\documentclass{amsart}
\usepackage{amssymb,latexsym}
\input epsf
\newtheorem{theorem}{Theorem}[section]
\newtheorem{proposition}[theorem]{Proposition}
\newtheorem{corollary}[theorem]{Corollary}

\newtheorem{lemma}[theorem]{Lemma}

\newcommand{\nc}{{\rm {\bf L}}}
\newcommand{\lk}{{\rm lk}}

\renewcommand{\to}{\rightarrow}

\newcommand{\sm}{{\setminus}}
\begin{document}
\title[Enumerative aspects of generalized associahedra]
{On some enumerative aspects of generalized associahedra}

\author{Christos~A.~Athanasiadis}
\address{Department of Mathematics\\
University of Crete\\
71409 Heraklion, Crete, Greece}
\email{caa@math.uoc.gr}

\date{August 1, 2005.}
\begin{abstract}
We prove a conjecture of F. Chapoton relating certain enumerative
invariants of (a) the cluster complex associated by S. Fomin and 
A. Zelevinsky to a finite root system and (b) the lattice of 
noncrossing partitions associated to the corresponding finite 
real reflection group. 
\end{abstract}

\maketitle

\section{The result}
\label{intro}

Let $\Phi$ be a finite root system spanning an $n$-dimensional Euclidean 
space $V$ with corresponding finite reflection group $W$. Let $\Phi^+$ be 
a positive system for $\Phi$ with corresponding simple system $\Pi$. The 
cluster complex $\Delta (\Phi)$ was introduced by S. Fomin and A. 
Zelevinsky within the context of their theory of cluster algebras 
\cite{FZ1, FZ2, FZ3}. It is a pure $(n-1)$-dimensional simplicial 
complex on the vertex set $\Phi^+ \cup (-\Pi)$ which is homeomorphic to 
a sphere \cite{FZ2}. It was realized explicitly in \cite{CFZ} as the 
boundary 
complex of an $n$-dimensional simplicial convex polytope $P(\Phi)$, known 
as the simplicial generalized associahedron associated to $\Phi$. Although 
$\Delta (\Phi)$ was initially defined under the assumption that $\Phi$ is 
crystallographic \cite{FZ2}, its definition and main properties are valid 
without this restriction \cite[Section 5.3]{FR} \cite{FR2}; see \cite{FR} 
for an expository treatment of generalized associahedra.

The combinatorics of $\Delta (\Phi)$ is closely related to that of a 
finite poset $\nc_W$, known as the lattice of noncrossing partitions  
associated to $W$ \cite{Be, BW1} (see Section \ref{pre} for 
definitions). It is known, for instance (see \cite[Theorem 5.9]{FR}), 
that the $h$-polynomial of $\Delta (\Phi)$ is equal to the rank 
generating polynomial of $\nc_W$. In particular the number of facets of 
$\Delta (\Phi)$ is equal to the cardinality of $\nc_W$. This number is 
a Catalan number if $\Phi$ has type $A_n$ in the Cartan-Killing 
classification; in that case $P(\Phi)$ is the polar polytope to the 
classical $n$-dimensional associahedron \cite[Section 3.1]{FR} and 
$\nc_W$ is 
isomorphic to the lattice of noncrossing partitions of the set $\{1, 
2,\dots,n+1\}$ \cite[Section 5.1]{FR}. The poset $\nc_W$ is a self-dual 
graded lattice of rank $n$ which plays an important role in the geometric 
group theory and topology of finite-type Artin groups; see \cite{Mc} for
a related survey article. 

The $F$-triangle for $\Phi$, introduced by F. Chapoton \cite[Section 
2]{Ch1}, is a refinement of the $f$-vector of $\Delta (\Phi)$ defined 
by the generating function
\begin{equation}
F(\Phi) = F(x, y) \ = \ \sum_{k=0}^n \sum_{\ell=0}^n \ f_{k, \ell} 
\, x^k y^\ell
\label{eq:F}
\end{equation}
where $f_{k, \ell}$ is the number of faces of $\Delta (\Phi)$ 
consisting of $k$ positive roots and $\ell$ negative simple roots.
Clearly $f_{k, \ell} = 0$ unless $k + \ell \le n$. The $M$-triangle 
for $W$ is defined similarly \cite[Section 3]{Ch1} as 
\begin{equation}
M(W) = M(x, y) \ = \ \sum_{a \preceq  b} \ \mu (a,b) \, x^{r(b)} 
y^{r(a)}
\label{eq:M}
\end{equation}
where $\preceq$ denotes the order relation of $\nc_W$, $\mu$ stands 
for its M\"obius function \cite[Section 3.6]{Sta}, $r(a)$ is the rank 
of $a \in \nc_W$ and the sum runs over all pairs $(a, b)$ of elements 
of $\nc_W$ with $a \preceq b$. The main objective of this note is to 
prove the following theorem, the rather surprising statement of which  
appears as \cite[Conjecture 1]{Ch1} and includes many of the known 
similarities between the enumerative properties of $\Delta (\Phi)$ 
and $\nc_W$ as special cases; see \cite[Sections 3.1--3.5]{Ch1}.
\begin{theorem} 
The $F$-triangle for $\Phi$ and $M$-triangle for $W$ are related by
the equality
\begin{equation}
(1-y)^n \, F(\frac{x+y}{1-y}, \frac{y}{1-y}) \ = \ M(-x,-y/x). 
\label{eq:main}
\end{equation}
\label{thm:main}
\end{theorem}

The proof of Theorem \ref{thm:main} (Section \ref{proof}) relies
on two known special cases, one relating the number $f_{n, 0}$ of facets 
of $\Delta (\Phi)$ consisting of only positive roots to the M\"obius 
number of $\nc_W$ \cite[(23)]{Ch1} and the one already mentioned, relating 
the $h$-polynomial of $\Delta (\Phi)$ to the rank generating polynomial 
of $\nc_W$. A case-free proof of the relevant statement in 
the former case was given in \cite[Corollary 4.4]{ABW} but no such 
proof exists at present in the latter case (see Remark 9.4 in \cite{Re} 
for an interesting attempt). To extract the proof of the theorem from 
the two special cases we utilize the appearance of the cluster complex 
in the context of the lattice $\nc_W$ in the work of T.~Brady and 
C.~Watt \cite{BW2}. This connection is briefly outlined in Section 
\ref{pre}, where the two special cases are conveniently generalized 
(Lemmas \ref{lem:mu(w)} and \ref{thm:links}). For root systems of type
$A$ and $B$, Theorem \ref{thm:main} follows from the computations in 
Sections 4 and 5 of \cite{Ch1} and results of C.~Krattenthaler \cite{Kr}.

\section{Noncrossing partitions and cluster complexes}
\label{pre}

Throughout this section $W$ is a finite real reflection group of rank 
$n$ with set of reflections $T$ and $\Phi$ is a root system spanning 
an $n$-dimensional Euclidean space $V$ with associated reflection 
group $W$. We refer the reader to \cite{Hu} and \cite{Sta} for 
background and any undefined terminology on root systems, finite 
reflection groups and partially ordered sets. 

\subsection{The lattice $\nc_W$.}
For $w \in W$ let $r(w) = r_T (w)$ denote the smallest integer $k$ such 
that $w$ can be written as a product of $k$ reflections in $T$. Define 
a partial order $\preceq$ on $W$ by letting 
\[ u \preceq v \ \ \ \text{if and only if} \ \ \ r (u) + r (u^{-1} 
v) = r (v), \]
in other words if there exists a shortest factorization of $u$ into 
reflections in $T$ which is a prefix of such a shortest factorization 
of $v$. Since $T$ is invariant under conjugation the function $r_T$ 
is constant on conjugacy classes of $W$ and we have $u \preceq v$ if and 
only if $wuw^{-1} \preceq wvw^{-1}$ for $u, v, w \in W$. 
\begin{lemma} 
Let $a, b, w$ be elements of $W$.
\begin{enumerate}
\itemsep=0pt
\item[(i)] $a \preceq aw \preceq b$ if and only if $w \preceq a^{-1} b 
\preceq b$.
\item[(ii)] $a \preceq aw \preceq b$ if and only if $a \preceq b w^{-1} 
\preceq b$.
\item[(iii)] $a \preceq b$ if and only if $a^{-1} b \preceq b$ and, in
that case, the interval $[a, b]$ is isomorphic to $[1, a^{-1} b]$. \qed
\end{enumerate}
\label{lem:nc}
\end{lemma}
\begin{proof}
Consider part (ii) and suppose that $a \preceq aw \preceq b$. From the 
definition of $\preceq$ we have $r(a) + r(w) = r(aw)$ and $b = awc$ with 
$r(aw) + r(c) = r(b)$. Let $c' = wcw^{-1}$ and observe that $r(c) = r(c')$. 
The factorization $b = ac'w$ implies that $r(ac') = r(a) + r(c')$, since
$r(ac') \le r(a) + r(c')$ on the one hand and $r(ac') \ge r(b) - r(w) =
r(a) + r(c) = r(a) + r(c')$ on the other. It follows that $a \preceq ac' 
= b w^{-1}$ and $b w^{-1} = ac' \preceq b$, so that $a \preceq b w^{-1} 
\preceq b$. The converse, as well as part (i), are treated in a similar 
way. Part (iii) follows from part (i). 
\end{proof}

The order $\preceq$ turns $W$ into a graded poset having the identity $1$ 
as its unique minimal element and rank function $r_T$. For $w \in W$ we 
denote by $\nc_W (w)$ the interval $[1, w]$ in this order. We are 
primarily interested in the case that $w$ is a Coxeter element $\gamma$ 
of $W$. Since all Coxeter elements of $W$ are conjugate to each other, the 
isomorphism type of the poset $\nc_W (\gamma)$ is independent of $\gamma$. 
This poset is denoted by $\nc_W$ when the choice of 
$\gamma$ is irrelevant and called the \emph{noncrossing partition lattice} 
associated to $W$. If $W$ is reducible, decomposing as a direct product
$W_1 \times W_2 \times \cdots \times W_k$ of irreducible parabolic 
subgroups, then $\nc_W$ is isomorphic to the direct product of the posets 
$\nc_{W_i}$.

\subsection{The cluster complex.}
\label{pre2}

It was shown in \cite[Section 8]{BW2} that the cluster complex $\Delta 
(\Phi)$ arises naturally in the context of the lattice $\nc_W$. We give 
a brief account of this connection here and refer the reader to \cite{BW2} 
for more details. We denote by $t_\alpha$ the reflection in the hyperplane 
in $V$ orthogonal to $\alpha$ and by $N$ the number of reflections in $T$. 
Let $\Phi^+$ be a fixed positive system for $\Phi$ and $\gamma$ be a 
corresponding bipartite Coxeter element of $W$, so that $\gamma = \gamma_+
\gamma_-$ where
\[ \gamma_\pm = \prod_{\alpha \in \Pi_\pm} t_\alpha \]
and $\Pi_+ = \{\alpha_1,\dots,\alpha_s\}$, $\Pi_- = 
\{\alpha_{s+1},\dots,\alpha_n\}$ are orthogonal sets which form a partition
of the simple system $\Pi$ determined by $\Phi^+$. Assume first that $\Phi$
is irreducible. Letting $\rho_i = t_{\alpha_1} t_{\alpha_2} \cdots 
t_{\alpha_{i-1}} (\alpha_i)$ for $i \ge 1$ (so that $\rho_1 = \alpha_1$), 
where the $\alpha_i$ are indexed cyclically modulo $n$, and $\rho_{-i} = 
\rho_{2N-i}$ for $i \ge 0$ we have 
\[ \begin{tabular}{ll}
$\{ \rho_1, \rho_2,\dots,\rho_N \} = \Phi^+$, \\ 
$\{ \rho_{N+i}: 1 \le i \le s\} = \{-\rho_1,\dots,-\rho_s \} = 
- \Pi_+$, \\ 
$ \{\rho_{-i}: 0 \le i < n-s \} = \{-\rho_{N-i}: 0 \le i < n-s\} =
- \Pi_-$.
\end{tabular} \]
Define an abstract simplicial complex $\Delta (\gamma)$ on the vertex 
set $\Phi_{\ge -1} = \Phi^+ \cup (-\Pi) = \{\rho_{-n+s+1},\dots,\rho_0, 
\rho_1,\dots,\rho_{N+s}\}$ by declaring a set $\sigma = \{\rho_{i_1}, 
\rho_{i_2},\dots,\rho_{i_k}\}$ with $i_1 < i_2 < \cdots < i_k$ to be 
a face if and only if 
\begin{equation}
w_\sigma = t_{\rho_{i_k}} t_{\rho_{i_{k-1}}} \cdots \ t_{\rho_{i_1}}
\label{eq:w_sigma}
\end{equation}
is an element
of $\nc_W (\gamma)$ of rank $k$. If $\Phi$ is reducible with irreducible 
components $\Phi_1, \Phi_2,\dots,\Phi_m$ then $\gamma = \gamma_1 \gamma_2 
\cdots \gamma_m$ where $\gamma_i$ is a bipartite Coxeter element for the 
reflection group $W_i$ corresponding to $\Phi_i$ and $\nc_W (\gamma)$ is
isomorphic to the direct product of the posets $\nc_{W_i} (\gamma_i)$. We 
define $\Delta 
(\gamma)$ as the simplicial join of the complexes $\Delta (\gamma_i)$, so 
that $\sigma \in \Delta (\gamma)$ if and only if $\sigma = \sigma_1 \cup 
\sigma_2 \cup \cdots \cup \sigma_m$ with $\sigma_i \in \Delta (\gamma_i)$
for $1 \le i \le m$. In that case we also define  
\begin{equation}
w_\sigma = w_{\sigma_1} w_{\sigma_2} \cdots w_{\sigma_m},
\label{eq2:w_sigma}
\end{equation}
where the $w_{\sigma_i}$ mutually commute, so that $w_\sigma$ is an element
of $\nc_W (\gamma)$ of rank equal to the cardinality of $\sigma$. 

Let $\Delta_+ (\gamma)$ and $\Delta_+ (\Phi)$ denote the induced 
subcomplexes of $\Delta (\gamma)$ and $\Delta (\Phi)$, respectively, on 
the vertex set $\Phi^+$. The following theorem is proved in \cite[Section 
8]{BW2} (see also \cite[Note 4.2]{BW2}) in the case of irreducible root
systems and extends by definition to the general case. 
\begin{theorem} {\rm (\cite{BW2})}
As an abstract simplicial complex $\Delta (\Phi)$ coincides with 
$\Delta (\gamma)$. In particular, $\Delta_+ (\Phi)$ coincides with 
$\Delta_+ (\gamma)$. \qed
\label{thm:BW}
\end{theorem}

As a consequence the complexes $\Delta (\gamma)$ and $\Delta_+ (\gamma)$
depend only on our fixed choice of $\Phi^+$. Recall (see \cite[Proposition 
1.6.4]{Be}) that any two reflections $t_1, t_2 \in T$ for which $t_1 t_2 
\preceq \gamma$ and $t_2 t_1 \preceq \gamma$ are commuting. This implies, 
in the case of an irreducible root system $\Phi$, that any rearrangement of 
the product in the right hand-side of (\ref{eq:w_sigma}) which is an 
element of $\nc_W (\gamma)$ of rank $k$ is equal to $w_\sigma$. In particular 
the map $\Delta (\gamma) \mapsto \nc_W$ sending $\sigma$ to $w_\sigma$ 
depends only on $\Phi^+$ and $\gamma$ and not on the specific linear orderings 
of $\Pi_+$ and $\Pi_-$ used in the definition of $\Delta (\gamma)$. Similar
remarks hold for the product in (\ref{eq2:w_sigma}) when $\Phi$ is reducible.
For any $w \in \nc_W (\gamma)$ the faces 
$\sigma$ of $\Delta_+ (\gamma)$ with $w_\sigma = w$ are the facets of a 
subcomplex $\Delta_+ (w)$ of $\Delta_+ (\gamma)$ (this corresponds to the 
complex denoted by $X(w)$ in \cite[Section 5]{BW2}). Clearly the sets of 
facets of the subcomplexes $\Delta_+ (w)$ for $w \in \nc_W (\gamma)$ form a 
partition of the set of faces of $\Delta_+ (\gamma)$ into mutually disjoint 
subsets.  

In the next proposition we gather some facts from \cite[Section 3]{FZ2},
modified according to some observations made in \cite[Section 8]{BW2}
(for the non-crystallographic root systems see, for instance, \cite{FR2}). 
For the sake of simplicity we assume that $\Phi$ is irreducible and let $R: 
\Phi_{\ge -1} \to \Phi_{\ge -1}$ be the map defined by
\[ R (\alpha) = \begin{cases}
\gamma^{-1} (\alpha), & \text{if $\alpha \notin \Pi_+ \cup (-\Pi_-)$}\\
-\alpha, & \text{if $\alpha \in \Pi_+ \cup (-\Pi_-)$}
\end{cases} \]
and let $R(\sigma) = \{R(\alpha): \alpha \in \sigma\}$ for $\sigma 
\subseteq \Phi_{\ge -1}$. For $\sigma \subseteq \Pi$ we denote by 
$\Phi_\sigma$ the standard parabolic root subsystem obtained by 
intersecting $\Phi$ with the linear span of $\Pi \sm \, \sigma$,
endowed with the induced positive system $\Phi_\sigma^+ = \Phi^+ \cap 
\Phi_\sigma$, and abbreviate $\Phi_\sigma$ as $\Phi_\alpha$ when 
$\sigma = \{\alpha\}$. 
\begin{proposition}
Let $\Phi$ be irreducible, $\alpha \in \Pi$ and $\sigma \subseteq 
\Phi_{\ge -1}$. 
\begin{enumerate}
\itemsep=0pt
\item[(i)] For $\sigma \in \Delta (\Phi)$ we have $-\alpha \in \sigma$ if 
and only if $\sigma \sm \, \{-\alpha\} \in \Delta (\Phi_\alpha)$. 
\item[(ii)] For any $\beta \in \Phi^+$ there exists $i$ such that $R^i 
            (\beta) \in (-\Pi)$.
\item[(iii)] $\sigma \in \Delta (\Phi)$ if and only if $R(\sigma) \in 
             \Delta (\Phi)$. 
\qed
\end{enumerate}
\label{prop:R}
\end{proposition}

\subsection{The M\"obius function.}
We write $\mu (w)$ instead of $\mu (\hat{0}, w)$ for the M\"obius function
between $\hat{0}$ and $w \in P$ of a finite poset $P$ with a unique minimal 
element $\hat{0}$. Let $\gamma$ be a bipartite Coxeter element of $W$, as 
in Section \ref{pre2}. It is known \cite[(23)]{Ch1} \cite[Corollary 
4.4]{ABW} that the M\"obius number $\mu (\gamma)$ of $\nc_W = \nc_W 
(\gamma)$ is equal, up to the sign $(-1)^n$, to the number of facets of 
$\Delta_+ (\Phi)$. This fact generalizes as follows. 
\begin{lemma} 
For $w \in \nc_W (\gamma)$ the number $(-1)^{r(w)} \mu (w)$ is equal to
the number of facets of $\Delta_+ (w)$.
\label{lem:mu(w)}
\end{lemma}
\begin{proof}
It suffices to treat the case that $\Phi$ is irreducible. By \cite[Corollary 
4.3]{ABW} $(-1)^{r(w)} \mu (w)$ is equal to
the number of factorizations $w = t_{\rho_{i_k}} t_{\rho_{i_{k-1}}} \cdots 
\ t_{\rho_{i_1}}$ of length $k = r(w)$ with $1 \le i_1 < i_2 < \cdots < i_k
\le N$. The set of such factorizations is in bijection with the set of faces
$\sigma$ of $\Delta_+ (\gamma)$ with $w_\sigma = w$, in other words with
the set of facets of $\Delta_+ (w)$.  
\end{proof}
The next corollary is the specialization $y=0$ of Theorem \ref{thm:main}.
\begin{corollary} 
The coefficient of $q^k$ in the characteristic polynomial 
\[ \chi (\nc_W, q) \ = \sum_{w \in \nc_W} \mu (w) \, q^{r(w)} \]
of $\nc_W$ is equal to $(-1)^k$ times the number of 
faces of $\Delta_+ (\Phi)$ of dimension $k-1$.
\label{cor:chi}
\end{corollary}
\begin{proof}
This follows from Theorem \ref{thm:BW}, Lemma \ref{lem:mu(w)} and the fact 
that the set of faces of $\Delta_+ (\gamma)$ of dimension $k-1$ is the 
disjoint union of the sets of facets of the subcomplexes $\Delta_+ (w)$ 
where $w$ ranges over all elements of $\nc_W (\gamma)$ of rank $k$.
\end{proof}

\subsection{Links and $h$-polynomials.}
Recall that the $h$-polynomial of an abstract simplicial complex $\Delta$
of dimension $n-1$ is defined as
\[ h(\Delta, y) \ = \ \sum_{i=0}^n \ f_i (\Delta) \, y^i (1-y)^{n-i} \]
where $f_i (\Delta)$ is the number of faces of $\Delta$ of dimension $i-1$.
The link of a face $\sigma$ of $\Delta$ is the abstract simplicial complex 
$\lk_\Delta (\sigma) = \{\tau \sm \, \sigma: \sigma \subseteq \tau \in 
\Delta\}$. It is known \cite[Section 3.1]{Ch1} \cite[Theorem 5.9]{FR} that 
\begin{equation}
h(\Delta (\Phi), y) \ = \sum_{a \in \nc_W} y^{r(a)}  
\label{eq:nar}
\end{equation}
for any root system $\Phi$ (the coefficients of either hand-side of 
(\ref{eq:nar}) are known as the \emph{Narayana numbers} associated to $W$). 
This fact generalizes as follows, where $\gamma$ is as in Section \ref{pre2}.
\begin{lemma} 
For any face $\sigma$ of $\Delta = \Delta (\Phi)$ we have
\[ h(\lk_\Delta (\sigma), y) \ = 
\sum_{a \in \nc_W (\gamma w_\sigma^{-1})} y^{r(a)}. \]
\label{thm:links}
\end{lemma}
\begin{proof}
We will use induction on the rank of $\Phi$. The proposed equality follows
by induction if $\Phi$ is reducible, since both hand-sides are multiplicative
in $W$, and reduces to (\ref{eq:nar}) if $\sigma$ is empty. Henceforth we
assume that $\Phi$ is irreducible and let $\sigma = \{\rho_{i_1}, 
\rho_{i_2},\dots,\rho_{i_k}\}$, in the notation of Section \ref{pre2}, with 
$i_1 < i_2 < \cdots < i_k$ and $k \ge 1$. We distinguish three cases.

\vspace{0.1 in}
\noindent
{\bf Case 1.} $\sigma \cap (-\Pi_-) \neq \emptyset$. Then $\rho_{i_1} 
= - \alpha$ with $\alpha \in \Pi_-$ and hence, by Proposition \ref{prop:R} 
(i), we have $\lk_\Delta (\sigma) = \lk_{\Delta'} (\sigma')$ where 
$\Delta' = \Delta (\Phi_\alpha)$ and $\sigma' = \sigma \sm \, \{-\alpha\}$. 
The induction hypothesis implies that 
\begin{equation}
h(\lk_\Delta (\sigma), y) \ = 
\sum_{a \in \nc_W (\gamma' w_{\sigma'}^{-1})} y^{r(a)} 
\label{eq:step}
\end{equation}
where $\gamma' = \gamma t_\alpha$. Clearly $t_{\rho_{i_k}} \cdots \
t_{\rho_{i_2}} = w_\sigma t_\alpha$ is a rank $k-1$ element of $\nc_W 
(\gamma')$ and hence $w_{\sigma'} = w_\sigma t_\alpha$. Thus $\gamma' 
w_{\sigma'}^{-1} = \gamma w_\sigma^{-1}$ and the result follows from 
(\ref{eq:step}).

\vspace{0.1 in}
\noindent
{\bf Case 2.} $\sigma \cap (-\Pi_+) \neq \emptyset$. Then $\rho_{i_k} 
= - \alpha$ with $\alpha \in \Pi_+$ and (\ref{eq:step}) continues to 
hold with $\gamma' = t_\alpha \gamma$ and $w_{\sigma'} = t_\alpha 
w_\sigma$, where $\Delta'$ and $\sigma'$ have the same meaning as in 
the previous case. Since 
$\gamma' w_{\sigma'}^{-1} = t_\alpha (\gamma w_\sigma^{-1}) \, t_\alpha$ 
is conjugate to $\gamma w_\sigma^{-1}$, the poset $\nc_W (\gamma' 
w_{\sigma'}^{-1})$ is isomorphic to $\nc_W (\gamma w_\sigma^{-1})$ and 
the result follows again from (\ref{eq:step}). 

\vspace{0.1 in}
\noindent
{\bf Case 3.} $\sigma \cap (-\Pi) = \emptyset$. Let $\ell \ge 1$ be 
such that $\rho_{i_j} \in \Pi_+$ if and only if $j < \ell$ and let $\sigma' 
= \{\rho_{i_\ell},\dots,\rho_{i_k}\}$. Proposition \ref{prop:R} (iii)
implies that $\lk_\Delta (\sigma)$ is isomorphic to $\lk_\Delta 
(R(\sigma))$. Since $R(\alpha) = - \alpha \in (-\Pi_+)$ for 
$\alpha \in \sigma \sm \, \sigma' = \sigma \cap \Pi_+$, we have in turn
that $\lk_\Delta (R(\sigma)) \ = \ \lk_{\Delta'} (R(\sigma'))$ by part (i) 
of the same proposition, where $\Delta' = \Delta (\Phi_{\sigma \sm \, 
\sigma'})$. Suppose first that $\sigma \cap \Pi_+$ is nonempty, so that 
$\ell \ge 2$ and $\Phi_{\sigma \sm \, \sigma'}$ has smaller rank than 
$\Phi$. The previous observations and the induction hypothesis imply that 
\begin{equation}
h(\lk_\Delta (\sigma), y) \ = 
\sum_{a \in \nc_W (\gamma' w_{R(\sigma')}^{-1})} y^{r(a)}
\label{eq2:step}
\end{equation}
where $\gamma' = t_{\rho_{i_1}} \cdots t_{\rho_{i_{\ell-1}}} \gamma$. 
Since $w_\sigma = t_{\rho_{i_k}} \cdots t_{\rho_{i_1}}$ is an element of 
$\nc_W (\gamma)$ of rank $k$ we have that $t_{\rho_{i_k}} \cdots 
t_{\rho_{i_\ell}}$ is an element of $\nc_W (\gamma \, t_{\rho_{i_1}} \cdots 
t_{\rho_{i_{\ell-1}}})$ of rank $k-\ell+1$ and hence that $t_{R(\rho_{i_k})} 
\cdots t_{R(\rho_{i_\ell})} = \gamma^{-1} t_{\rho_{i_k}} \cdots 
t_{\rho_{i_\ell}} \gamma$ is an element of $\nc_W (\gamma')$ of rank 
$k-\ell+1$. Therefore $w_{R(\sigma')} = \gamma^{-1} t_{\rho_{i_k}} \cdots 
t_{\rho_{i_\ell}} \gamma$ and $\gamma' w_{R(\sigma')}^{-1} = w_\sigma^{-1} 
\gamma$ is conjugate to $\gamma w_\sigma^{-1}$, so that the result follows 
from (\ref{eq2:step}) as in the second case. Finally suppose that $\sigma 
\cap \Pi_+ = \emptyset$. By the previous argument it suffices to prove that
\[ h(\lk_\Delta (R(\sigma)), y) \ = 
\sum_{a \in \nc_W (\gamma w_{R(\sigma)}^{-1})} y^{r(a)}. \]
In view of Proposition \ref{prop:R} (ii), applying similarly $R$ 
sufficiently many times brings us back either to the previous situation 
or to one of the first two cases.
\end{proof}

\section{Proof of Theorem \ref{thm:main}}
\label{proof}

Throughout this section $\gamma$ is a bipartite Coxeter element of $W$, as 
in Section \ref{pre2}, and $|\sigma|$ denotes the cardinality of a finite 
set $\sigma$.
                    
\vspace{0.1 in}
\noindent
\emph{Proof of Theorem \ref{thm:main}.} To simplify notation let us write
$\Delta$, $\Delta_+$ and $\nc$ instead of $\Delta (\Phi)$, $\Delta_+ (\Phi)$ 
and $\nc_W = \nc_W (\gamma)$, respectively, and $a \preceq_\nc b$ instead of
$a \preceq b$ with $a, b \in \nc$. From (\ref{eq:F}) we have
\bigskip
\[ \begin{aligned}
(1-y)^n \, F(\frac{x+y}{1-y}, \frac{y}{1-y}) 
& \ = \ \sum_{k, \ell} \ f_{k, \ell} (x+y)^k y^\ell (1-y)^{n-k-\ell} \\
& \ = \ \sum_{k, \ell, i} \ f_{k, \ell} {k \choose i} 
             \, x^i y^{k+\ell-i} (1-y)^{n-k-\ell} \\
& \ = \ \sum_{\tau \in \Delta} \sum_{\substack{\sigma \in \Delta_+ \\ 
                \sigma \subseteq \tau}} x^{|\sigma|} y^{|\tau| - |\sigma|} 
                (1-y)^{n-|\tau|}\\
& \ = \ \sum_{\sigma \in \Delta_+} x^{|\sigma|} 
        \, \sum_{\substack{\tau \in \Delta \\ \sigma \subseteq \tau}} 
         y^{|\tau| - |\sigma|} 
(1-y)^{n-|\tau|} \\
& \ = \ \sum_{\sigma \in \Delta_+} x^{|\sigma|} 
    \sum_{\tau' \in \lk_\Delta (\sigma)} y^{|\tau'|} 
(1-y)^{n-|\sigma|-|\tau'|}
\end{aligned} \]
and hence that
\begin{equation}
(1-y)^n \, F(\frac{x+y}{1-y}, \frac{y}{1-y}) \ = \ \sum_{\sigma \in 
\Delta_+} x^{|\sigma|} \, h(\lk_\Delta (\sigma), y). 
\label{eq:proof1}
\end{equation}
Similarly, using (\ref{eq:M}) and observing from Lemma \ref{lem:nc} (iii)
that $\mu (a,b) = \mu (w)$, where $w = a^{-1} b$, we have
\bigskip
\[ \begin{aligned}
M(-x,-y/x)
& \ = \ \sum_{a \preceq_\nc b} \ \mu (a,b) \ (-x)^{r(b) - r(a)} y^{r(a)} \\
& \ = \ \sum_{a \preceq_\nc aw} \ \mu (w) \ 
(-x)^{r(w)} y^{r(a)} \\
& \ = \ \sum_{w \in \nc} \ (-x)^{r(w)} \mu (w) 
\sum_{a \preceq_\nc \gamma w^{-1}} y^{r(a)}.
\end{aligned} \]
Observe that we have used Lemma \ref{lem:nc} (ii) to conclude that for 
$a, w \in \nc$ we have $a \preceq aw \in \nc$ if and only if $a 
\preceq \gamma w^{-1} \in \nc$. From the last expression and Lemma 
\ref{lem:mu(w)} we have
\[ M(-x,-y/x) \ = \ \sum_{w \in \nc} \ x^{r(w)} \sum_{\substack{\sigma 
\in \Delta_+ \\ w_\sigma = w}} \sum_{a \preceq_\nc \gamma w^{-1}} 
y^{r(a)} \]
or, equivalently, 
\begin{equation}
M(-x,-y/x) \ = \ \sum_{\sigma \in \Delta_+} \ x^{|\sigma|} 
\sum_{a \preceq_\nc \gamma w_\sigma^{-1}} y^{r(a)}.
\label{eq:proof2}
\end{equation}
In view of Lemma \ref{thm:links}, the result follows by comparing 
(\ref{eq:proof1}) and (\ref{eq:proof2}). 
\qed

%


\begin{thebibliography}{99}
%
\bibitem{ABW}
C.A.~Athanasiadis, T.~Brady and C.~Watt,
\textit{Shellability of noncrossing partition lattices},
preprint, 2005, 10 pages ({\tt ArXiV preprint math.CO/0503007}).
%
\bibitem{Be}
D.~Bessis,
\textit{The dual braid monoid},
Ann. Sci. Ecole Norm. Sup. {\bf~36} (2003), 647--683.
%
\bibitem{BW1}
T.~Brady and C.~Watt,
\textit{K($\pi$, 1)'s for Artin groups of finite type},
in \textit{Proceedings of the Conference on Geometric and Combinatorial
group theory, Part I (Haifa 2000)},
Geom. Dedicata {\bf~94} (2002), 225--250.
%
\bibitem{BW2}
T.~Brady and C.~Watt,
\textit{Lattices in finite real reflection groups},
preprint, 2005, 29 pages, Trans. Amer. Math. Soc. (to appear).
%
\bibitem{Ch1}
F.~Chapoton,
\textit{Enumerative properties of generalized associahedra},
S\'emin. Loth. de Combinatoire {\bf~51} (2004), $\#$ B51b (electronic).
%
%
\bibitem{CFZ}
F.~Chapoton, S.~Fomin and A.V.~Zelevinsky,
\textit{Polytopal realizations of generalized associahedra},
Canad. Math. Bull. {\bf~45} (2002), 537--566.
%
\bibitem{FR}
S.~Fomin and N.~Reading,
\textit{Root systems and generalized associahedra},
in \emph{Geometric Combinatorics}, IAS/Park City Mathematics Series 
(to appear).
%
\bibitem{FR2}
S.~Fomin and N.~Reading,
\textit{Generalized cluster complexes and Coxeter combinatorics}, 
preprint, 2005, 40pp ({\tt ArXiV preprint math.CO/0505085}).
%
\bibitem{FZ1}
S.~Fomin and A.V.~Zelevinsky,
\textit{Cluster algebras I: Foundations},
J. Amer. Math. Soc. {\bf~15} (2002), 497--529.
%
\bibitem{FZ2}
S.~Fomin and A.V.~Zelevinsky,
\textit{$Y$-systems and generalized associahedra},
Ann. of Math. {\bf~158} (2003), 977--1018.
%
\bibitem{FZ3}
S.~Fomin and A.V.~Zelevinsky,
\textit{Cluster algebras II: finite type classification},
Invent. Math. {\bf~154} (2003), 63--121.
%
\bibitem{Hu}
J.E.~Humphreys, 
Reflection groups and Coxeter groups,
Cambridge Studies in Advanced Mathematics {\bf~29}, 
Cambridge University Press, Cambridge, England, 1990.
%
\bibitem{Kr}
C.~Krattenthaler, personal communication, 2005. 
%
\bibitem{Mc}
J.~McCammond, 
\textit{An introduction to Garside structures}, 
preprint, 2004, 28 pages.
%
\bibitem{Re}
N.~Reading,
\textit{Clusters, Coxeter sortable elements and noncrossing partitions}, 
preprint, 2005, 28 pages ({\tt ArXiV preprint math.CO/0507186}).
%
\bibitem{Sta}
R.P.~Stanley,
Enumerative Combinatorics, vol.~1, Wadsworth \& Brooks/Cole,
Pacific Grove, CA, 1986; second printing, Cambridge University Press,
Cambridge, 1997.
%

\end{thebibliography}
\end{document}